\documentclass[12pt]{article}
\usepackage{amsthm}
\title{A combinatorial analysis of the average time for open-address hash coding insertion}
\author{Vaughan R. Pratt}
\begin{document}
\date{December 1973, revised 1974}
\maketitle
\begin{abstract}
In analysing a well-known hash-coding method, Knuth gave an exact expression for the average number of rejections encountered by players of a variant of musical chairs. We study a variant more closely related to musical chairs itself and deduce the same expression by a purely combinatorial approach.
\end{abstract}

In an analysis of the average time to insert an item when using open-address hash-coding, Knuth \cite[p. 528-530]{knuth1973} reduced the problem to the following question about musical chairs. Given $m$ chairs arranged in a circle and numbered clockwise from $0$ to $m-1$, if in turn each of $n$ people arrives at a randomly selected chair (his \emph{initial} chair) and walks clockwise until he finds an empty chair (his \emph{final} chair), what is the average number of \emph{rejections} (chairs found occupied during the search) per player?

For $m \geq n$, Knuth showed that the average number of rejections is $$\frac{1}{2n}\left( \frac{n^{\underline{2}}}{m} + \frac{n^{\underline{3}}}{m^{2}} + \frac{n^{\underline{4}}}{m^{3}} + \ldots \right)$$ where $n^{\underline{i}} = {n\choose{i}}i!$, the number of ways of choosing from $n$ items an ordered subset of $i$ items. Knuth's analysis is very complex and gives no clue as to why the answer is so simple. Knuth appreciated this issue and asked for an explanation \cite{knuth1973a}.

In the following analysis we reduce the problem yet further, to a game even closer to musical chairs, in which the average number of rejections is the same as in Knuth's game. We establish a one-to-one correspondence between the set of all rejections encountered, over the whole sample space, and a set whose cardinality is easily established. In this way we eliminate \emph{all} of the involved algebra of Knuth's proof. Regrettably, despite the absence of algebra, the following combinatorial proof is not as simple as the author first envisaged when he undertook the search for it. A simpler proof would be a valuable contribution to this facet of the theory of hash coding.

In musical chairs, the players move around until the music stops, at which time they all rush for chairs. To change Knuth's problem to resemble more closely this situation we shall require all $n$ people to have arrived before letting them search for their chairs. In effect, the people are partitioned into $m$ labelled blocks (the chairs providing the labels), some of which may be empty. All blocks then travel around the chairs clockwise at the same speed, each losing one member to each vacant chair it passes until it becomes empty. This process is the \emph{seating process}. To resolve seating conflicts, we rank the players in advance. The reader should have little difficulty in verifying that this version of musical chairs involves the same average number of rejections per player as does Knuth's version.

For our analysis, we shall count the total number $R$ of rejections suffered by all players, summed over all possible \emph{samples} (assignments of players to initial chairs). The average is then $\frac{R}{nm^{n}}$ (we must assume the players are distinguishable), which implies that $R$ should be $$\frac{1}{2}(n^{\underline{2}}m^{n-1} + n^{\underline{3}}m^{n-2} + n^{\underline{4}}m^{n-3} + ...)$$.

We say that a sample $S$ of $n$ players \emph{matches} a \emph{sub-sample} $T$ of $k$ out of those $n$ players, for $k \leq n$, when every member of the $j$-th block of $T$ is a member of the $j$-th block of $S$, for $j = 0, 1,..., m - 1$. When a sub-sample $P$ of $k$ out of our $n$ players has a block of two people assigned to chairs $c$, $k-2$ blocks, of one person each, assigned to chairs ${c+1,c+2,...,c+k-2}$ (addition mod $m$), and no other blocks, we call $P$ a \emph{$k$-pattern}. Exactly $m^{n-k}$ samples match any given $k$-pattern, and there are $\frac{n^{\underline{k}}m}{2}$ possible $k$-patterns, taking rotations into account. Hence the total number of possible matches of samples with $k$-patterns, for $k \geq 2$, is $$\frac{1}{2}(n^{\underline{2}}m^{n-1} + n^{\underline{3}}m^{n-2} + n^{\underline{4}}m^{n-3} + ...),$$ which equals the total number of rejections given by Knuth's analysis. We shall give a one-to-one correspondence between rejections and matches, which immediately gives us an alternative derivation of this expression for the number of rejections.

Given a sample $S$ and a particular rejection, say of player $a$ at the seat containing player $z$, we shall construct a particular match between some sample $T$ and some pattern $P$. We shall form $T$ from $S$ merely by relabelling the blocks of $S$, that is, by permuting them with respect to the chairs. The permutation amounts to gathering together a certain subset of the blocks, called the \emph{distinguished} blocks. The block containing player $a$ (call it $b_1$) remains where it is, at chair $c$, and is the first distinguished block. The other distinguished blocks of $S$ (call them $b_2, b_3,...,b_k$ for some $k \geq 1$) are slid around the circle counter-clockwise by interchanging them with the non-distinguished blocks, until blocks $b_1, b_2,...,b_k$ occupy chairs $c, c+1, ..., c+k-1$ respectively. The undistinguished blocks occupy the remaining chairs. Note that, among themselves, the $k$ distinguished blocks retain their relative initial seating order, as do the $m-k$ undistinguished blocks.

We now specify the distinguished blocks. Block $b_1$ has already been specified. Block $b_{i+1}$ is specified by the following procedure.
\begin{enumerate}
\item If block $b_i$ contains player $z$, then we are done and $k = i$. We refer to player $z$ as $p_k$.
\item Otherwise, let $d_i$ be the last chair to which $b_i$ lost a member (call him $p_i$) before $b_i$ (possibly exhausted) encountered $z$'s final chair. Then $b_{i+1}$ is the block whose initial chair is $d_i + 1$.
\end{enumerate}

This procedure completes the description of the construction of $T$. Note that there is a question as to whether $d_i$ is always defined, which is dealt with now.

In the following we let $[a, b]$ denote the set of chairs $\{a, a+1, a+2, ..., b\}$, the addition being modulo $m$. We let $(a,b] = [a,b] - \{a\}$, $[a,b) = [a,b] - \{b\}$, and $(a,b) = [a,b] - \{a, b\}$. We let $\underline{b}_i$ denote the initial chair of $b_i$, and \underline{z} the final chair of $z$. A player \emph{sits} in a set $S$ of chairs when his final chair is in $S$. A block of players sits in $S$ when \emph{some} player in that block sits in $S$. A set $C$ of chairs sits in $S$ when \emph{some} block whose initial chair is in $C$ sits in $S$. The dual of ``sits'' is \emph{sits only}; a block sits only in $S$ when \emph{all} players in the block sit in $S$, and a set $C$ of chairs sits only in $S$ when \emph{all} blocks whose initial chairs are in $C$ sit only in $S$.

The $b_i$'s determined by the procedure satisfy the following lemmas and theorem.

\newtheorem{Lemma 1}{Lemma 1}
\begin{Lemma 1}
$[\underline{b}_1,\underline{b}_k)$ and $[\underline{b}_k, \underline{z}]$ are disjoint.
\end{Lemma 1}
\begin{proof}
If not, then ${\underline{b}_1} \in [\underline{b}_k, \underline{z}]$  But then $a$, which is in $b_1$, would not be rejected by $z$ when $a$ reached $\underline{z}$, a contradiction.
\end{proof}

\newtheorem{Lemma 2}{Lemma 2}
\begin{Lemma 2}
$[\underline{b}_1,\underline{b}_k)$ does not sit in $[\underline{b}_k, \underline{z}]$.
\end{Lemma 2}
\begin{proof}
Each chair in $[\underline{b}_k, \underline{z}]$ is visited by an un-exhausted $b_k$ (because it contains $z$) before being visited by any block initially in $[\underline{b}_1,\underline{b}_k)$, by Lemma 1.
\end{proof}

\newtheorem{Theorem 3}{Theorem 3}
\begin{Theorem 3}
\renewcommand{\theenumi}{\roman{enumi}}
\renewcommand{\labelenumi}{\theenumi}
\mbox{}
\begin{itemize}
\item[(i)] $b_i$ is non-empty for $1 \leq i \leq k$.
\item[(ii)] $d_i \in [\underline{b}_i,\underline{b}_k)$ for $1 \leq i < k$.
\item[(iii)] $[\underline{b}_1, d_i]$ does not sit in $(d_i, \underline{b}_k]$ for $1 \leq i < k$ during the seating process.
\item[(iv)] $b_{i+1} \in (\underline{b}_i,\underline{b}_k]$ for $1 \leq i < k$.
\end{itemize}
\end{Theorem 3}

\begin{proof}
We use induction on $i$.
\begin{enumerate}
\item When $i = 1$, $b_1$ contains a. When $i>1$, $\underline{b}_i \in (\underline{b}_{i-1}, \underline{b}_k)$ by (iv). Hence a is rejected at $b_i$ as it progresses towards $\underline{z}$. The player seated at $\underline{b}_i$ must have originated in $[\underline{b}_1,\underline{b}_i]$ for $a$ to be rejected. By (iii) he cannot have originated in $[\underline{b}_1,d_{i-1}]$, which leaves only $\underline{b}_i$. Hence $b_i$ is non-empty.
\item By (i), $d_i$ exists. By the procedure , $d_i \in [\underline{b}_i, \underline{z}]$. By Lemma 2 and (iv), $d_i \notin [\underline{b}_k, \underline{z}]$. Hence $d_i \in [\underline{b}_i, \underline{b}_k]$.
\item By the procedure, $b_i$ does not sit in $(d_i, \underline{b}_k]$. Also $(\underline{b}_i, d_i]$ sits only in $(\underline{b}_i, d_i)$, otherwise $b_i$ would not  sit in $d_i$. So $[\underline{b}_i, d_i]$ does not sit in $(d_i, \underline{b}_k]$. By (iii), $[\underline{b}_1, d_{i-1}]$ does not sit in $(d_{i-1}, \underline{b}_k]$, and a portion does not sit in $(d_i, \underline{b}_k]$ by (ii). (For convenience take $[\underline{b}_1, d_0]$ to be empty.) Hence $[\underline{b}_1, d_i]$ does not sit in $(d_i, \underline{b}_k]$.
\item This follows directly from (ii) and the procedure, which makes $\underline{b}_{i+1} = d_i + 1$.
\end{enumerate}
\end{proof}

\newtheorem{Corollary 4}{Corollary 4}
\begin{Corollary 4}
The procedure halts.
\end{Corollary 4}

\begin{proof}
By (iv) of Theorem 3, as $i$ increases $b_i$ gets closer to $b_k$.
\end{proof}

We have now established that the distinguished blocks $b_1, b_2, ..., b_k$ are well-defined by the procedure. Hence the sample $T$ can now be constructed, by grouping together the distinguished blocks of $S$ as described earlier.

We construct the $(k+1)$-pattern $P$ to be matched by $T$ thus. Player $a$ is assigned to chair $c$ ($b_1$'s initial chair in both $S$ and $T$). For each distinguished block $b_i$ we assign player $p_i$ (as defined in the procedure) to chair $c+i-1$. This completes the construction of $P$. It is trivial to verify that $T$ matches $P$.

We have thus far exhibited a map from rejections to matches. To see that the map is a bijection, we show that its inverse is totally and uniquely defined.

Given a sample $T$ that matches a $(k+1)$-pattern $P$, we show how to reconstruct the corresponding sample $S$. The pattern serves to identify the blocks $b_1, b_2, ..., b_k$. To form $S$, we merge these blocks with the remaining $m-k$ undistinguished blocks of $T$. We leave block $b_1$ where it is, in chair $c$. Let player $a$ (the one who suffers the rejection we are constructing) be the second-ranked of the two players in $P$ at chair $c$, the other of which we call $p_1$. We assign undistinguished (possibly empty) blocks, in the order they appear directly following $b_k$ in $T$, to chairs $c+1, c+2, ...$ until $b_1$ and $b_2$ are sufficiently far apart that in the seating process for S player $p_1$ will be seated before $b_1$ arrives at $b_2$'s initial chair. (Recall that for seating purposes we ranked the players in advance.) This fixes the position for $b_2$. We now insert further undistinguished blocks from T between $b_2$ and $b_3$, to seat player $p_2$ (the player in chair $c+1$ in pattern $P$), and so on until the position for $b_k$ is determined. Player $p_k$ then becomes player $z$. It should be clear that when $m \geq n$ this procedure will never require more than the available number of undistinguished blocks. This completes the demonstration of a one-to-one correspondence between rejections and matches.

In conclusion, we have characterized the rejection aspect of musical chairs in terms of a reasonably natural correspondence with ``matching chairs'', in which we are interested in matches between samples and patterns, rather than in rejections of players by occupied chairs. The correspondence is arrived at by the way of an ``un-merging'' of blocks which is specified by using the notion of ``last player seated before rejection.''

\bibliographystyle{plain}
\bibliography{pratt}
\end{document}